\documentstyle[11pt,leqno]{article}
\newtheorem{theorem}{{\sc Theorem}}
\newcommand{\bt}{\begin{theorem}}
\newcommand{\et}{\end{theorem}}
\newcommand{\newsection}[1]{\setcounter{equation}{0} \setcounter{theorem}{0}
\section{#1}}

\newcommand{\NI}{\noindent}
\newcommand{\bea}{\begin{eqnarray}}
\newcommand{\eea}{\end{eqnarray}}

\def \spec#1 {\mathop{#1}}

\def \b #1 {\bf #1}

\newcommand {\CC}{\centerline}

\newcommand{\clf}{{\cal F}}

\newcommand{\clm}{{\cal M}}

\newcommand{\ity}{\infty}
\newcommand{\raro}{\rightarrow}

\newcommand{\vsp}{\vskip 1em}

\newcommand{\be}{\begin{equation}}
\newcommand{\ee}{\end{equation}}
\newcommand{\ben}{\begin{eqnarray*}}
\newcommand{\een}{\end{eqnarray*}}

\begin{document}
\vspace{5mm}
\centerline{\bf{Evaluation of Bid and Ask Pricing for European Option Under Mixed Fractional}}
\centerline{\bf{Brownian Motion Environment with Superimposed Jumps}}
\vsp
\CC{B.L.S. PRAKASA RAO}
\CC{C R Rao Advanced Institute of Mathematics, Statistics and Computer Science}
\CC{Hyderabad 500046, India}
\vspace{2mm}
\NI{\bf Abstract:} We investigate the valuation of the bid and ask prices for European option under mixed fractional Brownian motion environment in the presence of superimposed jumps by an independent Poisson process..
\vsp
\NI{Mathematics Subject Classification:} 60G22, 91B28.
\vsp
\newsection{Introduction} 
Prakasa Rao (2015) studied European option pricing for processes driven by mixed fractional Brownian motion with superimposed jumps. Pricing geometric Asian power options under mixed fractional Brownian motion environment is investigated in Prakasa Rao (2016). In a recent paper, Li and Wang (2021) studied valuation of bid and ask prices for European option under mixed fractional Brownian motion environment. The classical option pricing theory is built on the law of one price ignoring the effect of market liquidity on bid-ask spread. Cherney and Madan (2009) and Madan and Cherney (2010) developed the theory of conic finance which replaces the law of one price by the law of two prices allowing market participants sell to the market at the bid price and buy from the market at the higher ask price. The difference between the bid and ask prices is called the bid-ask spread which is an indication of the market liquidity. The bid price of a cash flow $X$ is defined by its discounted distorted expectation and the ask price by the negative of the discounted distorted expectation of the cash flow $-X.$ Madan and Cherney (2010) proposed to model the liquidity of the market by a single market stress level $0\leq \gamma \leq 1$ with $\gamma=0$ indicating that the bid-ask spread is zero. They obtained closed-form expressions for bid and ask prices for European options. Madan and Schoutens (2016,2017) investigated models with conic option pricing. Sonono and Mashele (2016) studied estimation of bid-ask prices for options on LIBOR based instruments. Guillaume et al. (2019) discussed implied liquidity risk premia in option markets. Our aim is now to propose a geometric mixed fractional Brownian motion model for the stock price process with possible  jumps superimposed by an independent Poisson process and investigate the bid and ask  prices for the  European options for such a model as the market liquidity and hence the bid-ask spread might also depend on the jumps in the stock market price processes due to various unexpected happenings. Some special cases are studied in detail. For more on conic option pricing, see Madan and Schoutens (2016). 
\vsp
\newsection{Preliminaries} 
Let $(\Omega,\clf,P)$ be a probability space. A real-valued  function $f: [0,1]\raro [0,1]$  is called a {\it distortion Function} if it is nondecreasing with $f(0)=0$ and $f(1)=1.$ The function $foP$ defined by
$$foP(A)=f(P(A)), A \in \clf$$
is called a {\it distortion} of the probability measure $P$ with respect to the distortion function $f.$ Let $\hat f(x)= 1-f(1-x), x \in [0,1].$ For any $X\in L^\ity(\Omega,\clf,P),$ define
$$E_f(X)= \int_0^\ity(1-f(P(X\leq x))dx-\int_{-\ity}^0f(P(X\leq x))dx.$$
Note that
$$E_f(X)= \int_{-\ity}^{\ity}x\;df(F_X(x))$$
where $F_X(.)$ is the distribution function of the random variable $X$ under the probability measure $P.$ The quantity $E_f(X)$ is called the {\it distorted expectation} of the random variable $X.$ If $f(x)=x,0\leq x \leq 1,$ then $E_f(X)$ is the usual expectation of the random variable $X.$ 
\vsp
Let $\Phi(.)$ denote the standard normal distribution function. The function $f_\gamma(u)= \Phi(\Phi^{-1}(u)+\gamma)$ is called the WANG-transform.
\vsp
Following the conic finance theory developed by Cherney and Madan (2009) (cf. Li and Wang (2021)), we say that a risk $X$ is acceptable or marketed if $E_Q[X]\geq 0$ for every $Q\in \clm$ where $\clm$ is a convex set consisting of test measures under which the expected cash-flow is non-negative. Under a large set $\clm,$ the number of acceptable risks will form a small set as there will be more tests to be passed. A cash-flow $X$ with distribution function $F_X(.)$ is acceptable if its  distorted expectation relative to a class of distortion functions $f_\gamma$ is non-negative i.e, 
$$\int_{-\ity}^{\ity} x\;df_\gamma(F_X(x))\geq 0.$$ 
The acceptability index $\gamma$ quantifies the degree of distortion. For a given distortion function $f_\gamma,\gamma \geq 0,$ under which cash-flows are evaluated , the constant $\gamma$ can be interpreted s the market liquidity level. A liquidity level of $\gamma=0$ indicates that there is no distortion which is an indication of perfect liquidity and hence complete market. In conic finance theory, Madan and Cherney (2010) assumed that the market is a party willing to accept all stochastic cash-flows $X$ with an acceptability level $\gamma.$ The ask price of a claim $a_\gamma(X)$ is defined as the smallest price for which the cash-flow $X$ of selling the claim is acceptable at level $\gamma$ for the market. The bid price of a claim $b_\gamma(X)$ is defined by the highest price for which the cash-flow of buying the claim is acceptable at level $\gamma$ for the market.
\vsp

Let $X$ be the cash-flow generated by the claim at the maturity time $T$ and interest rate $r$ compounded continuously.Then the bid and ask prices of payoff $X$ are given by
\bea
b_\gamma(X)& =&  \sup\{b: E_{f_\gamma}[e^{-rT}X-b]\geq 0\}\\\nonumber
&=& e^{-rT}E_{f_\gamma}[X]
\eea
and
\bea
a_\gamma(X)& =&  \inf\{a: E_{f_\gamma}[a-e^{-rT}X]\geq 0\}\\\nonumber
&=& -e^{-rT}E_{f_\gamma}[-X]
\eea
\vsp
Suppose we are interested in the price of a stock as it evolves over time. Let $S(t)$ denote the price of the stock at time $t$. The process $\{S(t), t \geq 0\}$ is said to be a geometric Brownian motion with drift parameter $\mu$ and volatility parameter $\sigma$ if for all $y\geq 0$ and $t\geq 0$, the random variable $$S(t+y)/S(y)$$
is independent of all prices up to time $y$ and if the random variable
$$\log(S(t+y)/S(y))$$
is a Gaussian random variable with mean $\mu t$ and variance $t \sigma^2.$ Suppose the stock price process is a geometric Brownian motion. It is known that, given the initial price $S(0),$ the expected value of the price at time $t$ depends on the parameters $\mu$ and $\sigma^2$ and
$$E[S(t)]=S(0) e^{\mu t+\frac{1}{2}\sigma^2 t}.$$
It is also known that the sample paths of a geometric Brownian motion are continuous almost surely and hence the geometric Brownian motion is not suitable for modeling stock price process if there are likely to be jumps. Consider a call option having the strike price $K$ and expiration time $t.$ Under the assumption that the stock price process follows the geometric Brownian motion and the interest rate $r$ does not change over time,  Black-Scholes formula gives the unique  no arbitrage cost of the European call option (cf. Ross(2003)). For a detailed explanation of ``arbitrage" opportunity, see Ross (2003). It is now known that certain time series are long range dependent and it was thought that a process driven by a fractional Brownian motion, in particular,  a geometric fractional Brownian motion may be used as a model for modeling stock price process. A fractional Brownian motion process with Hurst parameter $H\in (0,1]$ is an almost surely continuous centered Gaussian process with
$$Cov(B_t^H,B_s^H)= \frac{1}{2}(|t|^{2H}+|s|^{2H}-|t-s|^{2H}), t,s \in R.$$
For properties of a fractional Brownian motion, see Prakasa Rao (2010). It is easy to see that this process reduces to the Wiener process or the Brownian motion if $H=\frac{1}{2}.$ Attempts to model stock price process using the fractional Brownian motion as the driving force were not successful as such a modeling leads to  an arbitrage opportunity under the model which violates the fundamental assumption of mathematical finance modeling of no arbitrage opportunity or no free lunch (cf. Kuznetsov (1999)). It is known that the fractional Brownian motion $\{B_t^H, t \geq 0\}$ with Hurst index $H\in (0,1)$ is neither a Markov process nor a semimartingale except when $H=\frac{1}{2}.$ Cheridito (2000,2001) introduced the concept of mixed fractional Brownian motion for modeling stock price process. He showed that the sum of a Brownian motion and a non-trivial multiple of an independent fractional Brownian motion with Hurst index $H \in (0,1]$ is not a semimartingale for $H\in (0,\frac{1}{2})\cup (\frac{1}{2},\frac{3}{4}).$  However, if $H\in (\frac{3}{4},1],$  then the mixed fractional Brownian motion is equivalent to a multiple of Brownian motion and hence is  a semi-martingale. Hence, for $H\in (\frac{3}{4},1),$  the arbitrage opportunities can be excluded by modeling the stock price process $\{S_t,t \geq 0\}$ as the geometric mixed fractional Brownian motion given by
\be
S_t= S_0 \exp\{g(t)+ \sigma W_t^H+\epsilon W_t\}, t \geq 0
\ee
where $g(t)$ is a non-random function, $(\sigma, \epsilon) \neq (0,0),$ and the processes $\{W_t^H, t \geq 0\}$ and $\{W_t, t \geq 0\}$ are {\it independent} fractional Brownian motion and Brownian motion respectively. Cherdito (2000) showed that this model is arbitrage-free (cf. Mishura and Valkeila (2002)). It is known that the sample paths of the  mixed fractional Brownian motion or the geometric mixed fractional Brownian motion are almost surely continuous (cf. Zili (2006), Prakasa Rao (2010)).In order to take into account the long-memory property as well as to model the fluctuations in the stock prices in a financial market, one can use the mixed fractional Brownian motion as the driving force to model the stock price process. Sun (2013) discussed pricing currency options using  the mixed fractional Brownian motion model.
Sun and Yan (2012) discussed application of the mixed-fractional models to credit risk pricing. Yu and Lan (2008) studied the European call option pricing under the  mixed fractional Brownian environment. Since the sample paths of the geometric mixed fractional Brownian motion are almost surely continuous, it is not suitable for modeling stock price process with possible jumps. Mishura (2008) discussed sufficient conditions for the existence and uniqueness of solutions of stochastic differential equations driven by a mixed fractional Brownian motion. Yu and Yan (2008) derived  the analogue of the  Black-Scholes formula for the European call option price when the stock price process is the  geometric mixed fractional Brownian motion with interest rate $r$ is being constantly compounded continuously, the strike time is $t$ and the strike price is $K.$ Suppose the initial price of the stock at time 0 is $s.$ They showed that the European call option price is given by the formula
$$C(s,t,K, \sigma, r, \epsilon)= s\;\Phi(d_1)-Ke^{-rt}\;\Phi(d_2)$$
where
$$d_1= \frac{\log(s/K)+rt+\frac{1}{2}\sigma^2t^{2H}+\frac{1}{2}\epsilon^2t}{\sqrt{\sigma^2t^{2H}+\epsilon^2t}},$$
$$d_2= \frac{\log(s/K)+rt-\frac{1}{2}\sigma^2t^{2H}-\frac{1}{2}\epsilon^2t}{\sqrt{\sigma^2t^{2H}+\epsilon^2t}}.$$
and $\Phi(.)$ is the standard Gaussian distribution function.
\vspace {2mm}
\newsection{Adding Jumps to Geometric Mixed Fractional Brownian Motion}
\vsp
It is now known that modeling of the stock price process using the geometric Brownian motion is not useful as it does now allow possibility of a discontinuous price jump either in the upward or downward direction and to model long range dependence. Under the assumption of the geometric Brownian motion, the probability of having a jump is zero. Since such jumps do occur in practice for various reasons, it is important to consider a model for the stock price process that allows possibility of jumps in the process. 
\vsp
We assume that there are no transaction costs, trading is continuous and  the interest rate is constant and compounded continuously. We have indicated that there are no arbitrage opportunities under the mixed fractional Brownian motion model whenever the Hurst index $H\in (\frac{3}{4},1].$
\vsp
We would like to investigate the valuation of the bid and ask prices for European options under the mixed fractional Brownian motion with Hurst index $H>3/4$ capturing the underlying asset returns in real markets in the presence of jumps. We now propose a jump mixed fractional Brownian motion model to capture jumps or discontinuities, fluctuations in the stock price process and to take into account the long range dependence of the stock price process and obtain the  European call option price for such models.
\vsp
Mixed fractional Brownian motion with superimposed jumps can be used for pricing currency options (cf. Xiao et al. (2010)). It  is based on the assumption that the exchange rate returns are generated by a two-part stochastic process, the first part dealing with small continuous price movements  generated by a mixed fractional Brownian motion and the second part by large infrequent price jumps generated by a Poisson process. As has been pointed out by Foad and Adem (2014), modeling by this two-part process is in tune with the market in which major information arrives infrequently and randomly. In addition, this process provides a model through heavy tailed distributions for modeling empirically observed distributions of exchange rate changes.
\vsp
We now introduce the Poisson process as a model for the jump times in the stock price process. Let $N(0)=0$ and let $N(t)$ denote the number of jumps in the process that occur by time $t$ for $t >0.$ Suppose that  the process $\{N(t), t \geq 0\}$ is a Poisson process with  stationary independent increments. Under such a process, the probability that there is a jump in a time interval of length $h$ is approximately $\lambda h$ for $h$ small and the probability of more than one jump in a time interval of length $h$ is almost zero for $h$ sufficiently small. Furthermore the probability that there is a jump in an interval does not depend on the information about the earlier jumps. Suppose that when the $i$-th jump occurs, the price of the stock is multiplied by an amount $J_i$  and the random sequence $\{J_i, i \geq 1\}$ forms an independent and identically distributed  (i.i.d.) sequence of random variables. In addition, suppose that the random sequence $\{J_i, i \geq 1\}$ is independent of the times at which the jumps occur. Let $S(t)$ denote the stock price at time $t$ for $t \geq 0.$ Then
 $$S(t)= \tilde S(t) \Pi_{i=1}^{N(t)} J_i, t \geq 0$$
 where $\{\tilde S(t), t \geq 0\}$ is the geometric mixed fractional Brownian motion  modeled according to equation (2.3) specified earlier. Note that, if there is a jump in the price process at time $t$, then the jump is of size $J_i$ at the $i$-th jump. Let
 $$J(t)= \Pi_{i=1}^{N(t)}J_i, t \geq 0$$
and we define $ \Pi_{i=1}^{N(t)}J_i=1$ if $N(t)=0.$ Note that
the random variable $\log \frac{\tilde S(t)}{\tilde S(0)}$ has the Gaussian distribution with mean g(t) and variance $\epsilon^2 t + \sigma^2 t^{2H}.$
Note that $S(0)= \tilde S(0)$ is the initial stock price and we assume that it is non-random. Observe that
$$E[S(t)]=E[\tilde S(t)J(t)]=E[\tilde S(t)] E[J(t)]$$
by the independence of the random variables $\tilde S(t)$ and $J(t).$ Furthermore
\bea
E[\tilde S(t)]&= & \tilde S(0) E[\exp\{g(t)+\sigma W_t^H+\epsilon W_t\}]\\\nonumber
& = & \tilde S(0)\exp\{g(t)+ \frac{1}{2} \sigma^2 t^{2H}+\frac{1}{2} \epsilon^2 t\}\\\nonumber
\eea
by the independence of the processes $\{W_t^H, t \geq0 \}$ and $\{W_t, t \geq 0\}$ and the properties of Gaussian random variables. It is easy to check that
$$E[J(t)]= e^{-\lambda t (1-E[J_1])}$$
and
$$Var [J(t)]= e^{-\lambda t (1-E[J^2_1])} - e^{-2\lambda t (1-E[J_1])}.$$
In particular, the equations given above show that
$$E[S(t)]= \tilde S(0)e^{g(t)+ \frac{1}{2} \sigma^2 t^{2H}+\frac{1}{2} \epsilon^2 t}e^{-\lambda t (1-E[J_1])}.$$
Suppose the interest rate $r$ is compounded continuously. Then the future value of the stock price $S(0)$, after time $t$, should be $S(0)e^{rt}$ under any risk-neutral probability measure. Under the no arbitrage assumption, it follows that
$$\tilde S(0)\exp\{g(t)+ \frac{1}{2} \sigma^2 t^{2H}+\frac{1}{2} \epsilon^2 t -\lambda t (1-E[J_1]) \}=S(0)e^{rt}.$$
Since $S(0)=\tilde S(0),$ it follows that the stock price process should satisfy the relation
$$g(t)+ \frac{1}{2} \sigma^2 t^{2H}+\frac{1}{2} \epsilon^2 t-\lambda t (1-E[J_1])=rt$$
which implies that
$$g(t)= rt-\frac{1}{2} \sigma^2 t^{2H}-\frac{1}{2} \epsilon^2 t+\lambda t (1-E[J_1])$$
under the no arbitrage assumption. The price for an European call option with strike price $K,$ strike time $t,$ and interest rate $r$ compounded continuously is equal to
$$E[e^{-rt}(S(t)-K)_+]$$
where the expectation is computed with respect to the Gaussian distribution with mean
$$rt-\frac{1}{2} \sigma^2 t^{2H}-\frac{1}{2} \epsilon^2 t+\lambda t (1-E[J_1])$$
and variance
$$ \sigma^2t^{2H}+ \epsilon^2t.$$
Here $a_+=a$ if $a \geq 0$ and $a_+=0$ if $a<0.$ Let $R_t$ be a Gaussian random variable with mean $rt-\frac{1}{2} \sigma^2 t^{2H}-\frac{1}{2} \epsilon^2 t+\lambda t (1-E[J_1])$ and variance $\sigma^2t^{2H}+ \epsilon^2t.$ Note that the option price for an European call option under this model is
\bea
E[e^{-rt}(S(t)-K)_+]&= &e^{-rt}E[(J(t)\tilde S(t)-K)_+]\\\nonumber
&= & e^{-rt} E[(J(t)\tilde S(0)e^{R_t}-K)_+]\\\nonumber
\eea
where $\tilde S(0)=S(0)$ is the initial price of the stock. 
\vsp
Let us now consider a European call option $C=(S_T-K)^+$ and a put option $P=(K-S_T)^+$ with strike price $K$ and maturity $T$ on the underlying price $\{S_t, t\geq 0\}$ following the model $\{S_t=\tilde S(t)J(t), t \geq 0\}$
with 
$$d\tilde S(t)=r\tilde S(t)+\sigma \tilde S(t)dW_t+\epsilon \tilde S(t)dW_t^H, t \geq 0$$
where $r$ denotes the risk-free interest rate, $\{W_t, t \geq 0\}$ is the standard Brownian motion and $\{W^H_t, t \geq 0\}$ is an independent standard fractional Brownian motion and the process $\{J(t), t \geq 0\}$ is independent process as defined above. By using the distorted expectation defined earlier, the bid price for the European call option is given by
\bea
b_\gamma(C) &= &e^{-rT}E_{f_\gamma}[(S_T-K)^+]\\\nonumber
&=& e^{-rT}E_{f_\gamma}[(\tilde S(T)J(T)-K)^+]\\\nonumber
&=& e^{-rT}\int_K^\ity(x-K)df_\gamma(F_{\tilde S(T)J(T)}(x))\\\nonumber
&=& e^{-rT}\int_K^\ity x df_\gamma(F_{\tilde S(T)J(T)}(x))-e^{-rT}\int_K^\ity K df_\gamma(F_{\tilde S(T)J(T)}(x)).\\\nonumber
\eea
Observe that
\bea
F_{\tilde S(T)J(T)}(x)&=& P(\tilde S(T)J(T)\leq x)\\\nonumber
&=& P(S_0\exp(rT+\sigma W_T+\epsilon W_T^H-\frac{1}{2}\sigma^2T-\frac{1}{2}\epsilon^2 T^{2H})J(T)\leq x).\\\nonumber
\eea
By similar discussion, the ask price for the European call option is
\bea
\;\;\;\;\\\nonumber
a_\gamma(C) &= &-e^{-rT}E_{f_\gamma}[-(S_T-K)^+]\\\nonumber
&=& -e^{-rT}E_{f_\gamma}[(\tilde S(T)J(T)-K)^+]\\\nonumber
&=& -e^{-rT}\int_{-\ity}^0x\;df_\gamma(1-F_{\tilde S(T)J(T)}(K-x))\\\nonumber
&=& -e^{-rT}\int_0^{\ity}x\;df_\gamma(1-F_{\tilde S(T)J(T)}(K+x))\\\nonumber
&=& -e^{-rT}\int_K^{\ity}(x-K)\;df_\gamma(1-F_{\tilde S(T)J(T)}(x))\\\nonumber
&=& -e^{-rT}\int_K^{\ity}x\;df_\gamma(1-F_{\tilde S(T)J(T)}(x))+e^{-rT}\int_K^{\ity}Kdf_\gamma(1-F_{\tilde S(T)J(T)}(x)).\\\nonumber
\eea
Applying the arguments similar to those given above, it can be checked that the bid and ask prices for the European put option are given by
the following equations:
\bea
\;\;\;\;\\\nonumber
b_\gamma(P) &= &e^{-rT}E_{f_\gamma}[(K-S_T)^+]\\\nonumber
&=& e^{-rT}E_{f_\gamma}[(K-\tilde S(T)J(T))^+]\\\nonumber
&=& e^{-rT}\int_0^\ity x\; df_\gamma(1-F_{\tilde S(T)J(T)}(K-x))\\\nonumber
&=& -e^{-rT}\int_0^K (K-x)\; df_\gamma(1-F_{\tilde S(T)J(T)}(x))\\\nonumber
&=& -e^{-rT}\int_0^K K\; df_\gamma(1-F_{\tilde S(T)J(T)}(x))+e^{-rT}\int_0^K x \;df_\gamma(1-F_{\tilde S(T)J(T)}(x))\\\nonumber
\eea
and
\bea
a_\gamma(P) &= &-e^{-rT}E_{f_\gamma}[-(K-S_T)^+]\\\nonumber
&=& -e^{-rT}E_{f_\gamma}[-(K-\tilde S(T)J(T))^+]\\\nonumber
&=& -e^{-rT}\int_{-\ity}^0 x\; df_\gamma(F_{\tilde S(T)J(T)}(K+x))\\\nonumber
&=& -e^{-rT}\int_0^\ity x \;df_\gamma(F_{\tilde S(T)J(T)}(K-x))\\\nonumber
&=& e^{-rT}\int_0^K (K-x)\; df_\gamma(F_{\tilde S(T)J(T)}(x))\\\nonumber
&=& e^{-rT}\int_0^K K\; df_\gamma(F_{\tilde S(T)J(T)}(x))-e^{-rT}\int_0^K x\; df_\gamma(F_{\tilde S(T)J(T)}(x)).\\\nonumber
\eea
\vsp
The bid and ask prices for call and put options depend on the parameter $\gamma.$ In case $\gamma=0,$ then the bid price is equal to the ask price and the problem reduces to the one-price situation as discussed in Prakasa Rao (2015,2016) as bid-ask spread is zero.
\newsection{Special case}
Let us consider a special case of the model for the stock price discussed earlier where we can obtain the bid and ask prices in closed form expressions. Suppose that the jumps $\{J_i, i \geq 1\}$ are i.i.d. log-normally distributed with parameters $\mu_1$ and $\sigma_1^2.$ It is easy to see that
$$E[J_1]=  e^{\mu_1+\frac{1}{2}\sigma_1^2}.$$
Let $X_i= \log J_i, i \geq 1.$ Then the random variables $X_i, i \geq 1$ are i.i.d. with Gaussian distribution with mean $\mu_1$ and variance $\sigma_1^2.$ Observe that
$$J(t)= e^{\sum_{i=1}^{N(t)}X_i}$$
in this special case. The distribution function of $S_T$ is given by 
\bea
\;\;\;\;\\\nonumber
F_{S_T}(x)&=& P(S_T \leq x)\\\nonumber
&=& P(\tilde S(T)J(T)\leq x)\\\nonumber
&=& P(S_0\exp(rT+\sigma W_T+\epsilon W_T^H-\frac{1}{2}\sigma^2T-\frac{1}{2}\epsilon^2T^{2H}+\log J(T))\leq x)\\\nonumber
&=& P(S_0\exp(rT+\sigma W_T+\epsilon W_T^H-\frac{1}{2}\sigma^2T-\frac{1}{2}\epsilon^2T^{2H}+\sum_{i=1}^{N(T)}X_i)\leq x)\\\nonumber
&=& \sum_{n=0}^{\ity}P(S_0\exp(rT+\sigma W_T+\epsilon W_T^H-\frac{1}{2}\sigma^2T-\frac{1}{2}\epsilon^2T^{2H}+\sum_{i=1}^{N(T)=n}X_i)\leq x|N(T)=n)P(N(T)=n)\\\nonumber
&=& \sum_{n=0}^{\ity}P(S_0\exp(R_T+\sum_{i=1}^{N(T)=n}X_i)\leq x|N(T)=n)P(N(T)=n)\\\nonumber
\eea
where
\be
R_T= rT+\sigma W_T+\epsilon W_T^H-\frac{1}{2}\sigma^2T-\frac{1}{2}\epsilon^2T^{2H}.
\ee 
Since the Poisson process $\{N(t), t \geq 0\}$ is independent of the stock price process by assumption, it follows that
\bea 
F_{S_T}(x) &= &\sum_{n=0}^{\ity}P(S_0\exp(R_T+\sum_{i=1}^{n}X_i)\leq x)P(N(T)=n)\\\nonumber
&=& \sum_{n=0}^{\ity} P(S_0\exp(R_T+\sum_{i=1}^{n}X_i)\leq x)e^{-\lambda T}(\lambda T)^n (n!)^{-1}\\\nonumber 
&=& \sum_{n=0}^{\ity}P(R_T+\sum_{i=1}^{n}X_i \leq \log x -\log S_0)e^{-\lambda T}(\lambda T)^n (n!)^{-1}.\\\nonumber
\eea
Observe that the random variable $R_T+\sum_{i=1}^{n}X_i$ has the Gaussian distribution with the mean $rT-\frac{1}{2}\sigma^2 T-\frac{1}{2}\epsilon^2 T^{2H}+n \mu_1$ and the variance $\sigma^2T+\epsilon^2T^{2H}+n\sigma_1^2.$ Hence
\bea
\;\;\;\;\\\nonumber
F_{S_T}(x)&=& \sum_{n=0}^{\ity} P(Z\leq \frac{\log x-\log S_0-(rT-\frac{1}{2}\sigma^2T-\frac{1}{2}\epsilon^2T^{2H}+n\mu_1)}{\sqrt{T+T^{2H}+n\sigma_1^2}})e^{-\lambda T}(\lambda T)^n (n!)^{-1}
\eea
where $Z$ is a standard normal random variable.Hence
\bea
\;\;\;\;\\\nonumber
F_{S_T}(x)&=& \sum_{n=0}^{\ity}\Phi(\frac{\log x-\log S_0-(rT-\frac{1}{2}\sigma^2 T-\frac{1}{2}\epsilon^2 T^{2H}+n\mu_1)}{\sqrt{\sigma^2T+\epsilon^2T^{2H}+n\sigma_1^2}})e^{-\lambda T}(\lambda T)^n (n!)^{-1}\\\nonumber
&=& \sum_{n=0}^\ity G_n(x)e^{-\lambda T}(\lambda T)^n (n!)^{-1}\\\nonumber
\eea
where
\be
G_n(x)=\Phi(\frac{\log x-\log S_0-(rT-\frac{1}{2}\sigma^2T-\frac{1}{2}\epsilon^2T^{2H}+n\mu_1)}{\sqrt{\sigma^2T+\epsilon^2T^{2H}+n\sigma_1^2}}).
\ee
Note that the WANG-transform of the function $G_n(x)$ is
\bea
f_\gamma(G_n(x))&=&\Phi(\Phi^{-1}(G_n(x))+\gamma)\\\nonumber
&=& \Phi(\frac{\log x-\log S_0-(rT-\frac{1}{2}\sigma^2T-\frac{1}{2}\epsilon^2T^{2H}+n\mu_1)}{\sqrt{T+T^{2H}+n\sigma_1^2}}+\gamma).\\\nonumber
\eea
\vsp
\NI{\bf Acknowledgment:} This work was supported under the scheme "INSA Honorary Scientist " at the CR Rao Advanced Institute for Mathematics, Statistics and Computer Science, Hyderabad, India.
\vsp
\noindent{References;}

\begin{description}

\item Cheridito, P. (2000) {\it Regularizing fractional Brownian motion with a view toward stock price modeling}, Ph.D. Dissertation, ETH, Zurich.

\item Cheridito, P. (2001) Mixed fractional Brownian motion, {\it Bernoulli}, {\bf 7}, 913-934.

\item Cherney, A., and Madan, D.B. (2009) New measures of performance evaluation, {\it Review of Financial Studies}, {\bf 22}, 2571-2606.

\item Foad, S. and Adem, K. (2014) Pricing currency option in a mixed fractional Brownian motion with jumps environment, {\it Mathematical Problems in Engineering},  Vol. 2014, pp. 13, Article ID 858210, http://dx.doi.org/10.1155/2014/858210.

\item Guillaume, F., Junike, G., Leoni,P., Schoutens, W. (2019) Implied liquidity risk premia in option markets, {\it Annals of Finance}, {\bf 15}, 233-246. 

\item Kuznetsov, Yu. (1999) The absence of arbitrage in a model with fractal Brownian motion, {\it Russian Math. Surveys}, {\bf 54}, 847-848.

\item Li, Z. and Wang, X.-T. (2021) Valuation of bid and ask prices for European options under mixed fractional Brownian motion, Preprint.

\item Madan, D.B. and Cherney, A. (2010) Markets as a counterparty: an introduction to conic finance: an introduction to conic finance, {\it International Journal of Theoretical and Applied Finance}, {\bf 13}, 1149-1177.

\item Madan, D,B. and Schoutens, W. (2016) {\it Applied Conic Finance}, Cambridge University Press, Cambridge.

\item Madan, D,B. and Schoutens, W. (2017) Conic option pricing, {\it Journal of Derivatives}, {\bf 25}, 10-36.

\item Mishura , Y. (2008) {\it Stochastic Calculus for Fractional Brownian Motion and Related Processes}, Lecture Notes in Mathematics No. 1929, Springer, Berlin.

\item Mishura, Y.  and Valkeila, E. (2002) The absence of arbitrage in a mixed Brownian-fractional Brownian model, {\it Proceedings of the Steklov Institute of Mathematics}, {\bf 237}, 224-233.

\item Prakasa Rao, B.L.S. (2010) {\it Statistical Inference for Fractional Diffusion Processes}, Wiley, Chichester, UK.

\item Prakasa Rao, B.L.S. (2015) Option pricing for processes driven by mixed fractional Brownian motion with superimposed jumps, {\it Probability in the Engineering and Information Sciences}, {\bf 29}, 589-596.

\item Prakasa Rao, B.L.S. (2016) Pricing geometric Asian power options under mixed fractional Brownian motion environment, {\it Physica A: Statistical Machanics and its Applications}, {\bf 446}, 92-99.

\item Ross, Sheldon M. (2003) {\it An Elementary Introduction to Mathematical Finance,}, Second Edition, Cambridge University Press, Cambridge, UK.

\item Sonono, M.E. and Mashele, H.P. (2016) Estimation of bid-ask prices for options on LIBOR based instruments, {\it Finance Research Letters}, {\bf 19}, 33-41.

\item Sun, L. (2013) Pricing currency options in the mixed fractional Brownian motion, {\it Physica A, Vol. 392, no.16}, 3441-3458.

\item Sun, Xichao and Yan, Litan (2012) Mixed-fractional models to credit risk pricing, {\it Journal of Statistical and Econometric Methods}, {\bf 1}, 79-96.

\item Xiao, W.-L., Zhang, W,-G., Zhang, X.-L., and Wang, Y.-L. (2010). Pricing currency options in a fractional Brownian motion with jumps, {\it Economic Modelling}, {\bf 27}, 935-942.

\item Yu, Z. and Yan, L. (2008) European call option pricing under a mixed fractional Brownian motion environment, {\it Journal of University of Science and Technology of Suzhou}, {\bf 25}, 4-10 (in Chinese).

\item Zili, M. (2006) On the mixed fractional Brownian motion, {\it J. Appl. Math. Stochastic Anal.}, DOI10.1155/JAMSA/2006/32435.

\end{description}

\end{document}